\documentclass{amsart}
\usepackage{amssymb,amsfonts,latexsym}

\newtheorem{theorem}{Theorem}[section]
\newtheorem{lemma}[theorem]{Lemma}
\newtheorem{proposition}[theorem]{Proposition}
\newtheorem{corollary}[theorem]{Corollary}
\theoremstyle{definition}

\newtheorem{remark}[theorem]{Remark}

\newcommand{\Z}{{\mathbb Z}}

\numberwithin{equation}{section}

\begin{document}

\title[FORMULA FOR FINITE GROUP ACTION]{AN EXPLICIT FORMULA FOR THE ACTION OF A FINITE GROUP ON A COMMUTATIVE
RING}
\author{EHUD MEIR}
\address{Department of Mathematics, Technion-Israel Institute of
Technology, Haifa 32000, Israel}
\email{ehudm@tx.technion.ac.il}
\date{January 27, 2007}
\maketitle

\begin{abstract}
Let $G$ be a finite group, $k$ a commutative ring upon which $G$
acts. For every subgroup $H$ of $G$, the trace (or norm) map
$tr_H:k\rightarrow k^H$ is defined. $tr_H$ is onto if and only if
there exists an element $x_H$ such that $tr_H(x_H)=1$. We will show
that the existence of $x_P$ for every subgroup $P$ of prime order
determines the existence of $x_G$ by exhibiting an explicit formula
for $x_G$ in terms of the $x_P$, where $P$ varies over prime order
subgroups. Since $tr_P$ is onto if and only if $tr_{gPg^{-1}}$ is,
where $g\in G$ is an arbitrary element, we need to take only one $P$
from each conjugacy class. We will also show why a formula with less
factors does not exist, and show that the existence or non existence
of some of the $x_P$'s (where we consider only one $P$ from each
conjugacy class) does not affect the existence or non existence of
the others.
\end{abstract}

\begin{section}{INTRODUCTION}\label{introduction}
Let $G$ be a finite group, and let $k$ be a unital commutative ring
upon which $G$ acts as a group of automorphisms; i.e., we have a
group homomorphism $t:G\rightarrow \textrm{Aut}(k).$ For every
subgroup $H$ we define the trace map for this action: $$tr_H : k
\rightarrow k$$\begin{equation}tr_H(x) = \sum_{\sigma\in
H}{\sigma(x)}.\end{equation} For every $x\in k$, $tr_H(x)$ is an
$H$-invariant element of $k$, so we can consider the trace map as
$tr_H:k\rightarrow k^H$. The map $tr_H$ is a $k^H$-linear map, and
therefore it is onto if and only if there exists an element $x\in k$
such that $tr_H(x_H)=1$. The element $x_H$, if it exists, is not
unique in general. The reason why the surjectivity of the trace map
is of any interest to us is that its surjectivity is equivalent to
the projectivity of $k$ as a module over a certain skew group ring.
Details will be given in Section \ref{preliminaries}. From now on,
wherever we shall write $x_H$ we shall mean an element in $k$ which
satisfies the equation $tr_H(x_H)=1$. We will show later that if
$H<G$ and $tr_G$ is onto, then $tr_H$ is onto as well. In this paper
we will study the other direction. Namely, suppose we know that
$tr_H$ is onto for some of the subgroups $H$ of $G$. When can we
deduce that $tr_G$ is onto as well? It is known that if $H<G$ is any
subgroup, and $tr_K$ is onto for $K=H$ and for $K=N$ for every
subgroup $N$ such that $N \cap H=1$, then $tr_G$ is onto. See
\cite{Al1} for details. Moreover, In \cite{Al1} the following
formula for $x_G$ in terms of $x_H$ and $x_N$ is given:
\begin{equation}\label{formula1}x_G= \sum_{i=1}^s\Big[{x_{N_i}\Big(\prod_{t=1}^l{g_{it}(x_H)}\Big)}\Big].\end{equation}
The $N_i$ are subgroups of $G$ which intersects $H$ trivially, and
the $g_{it}$ are certain elements of $G$.

Using the above formula iteratively will give us a formula for $x_G$
in terms of the minimal subgroups (with respect to inclusion) of $G$.
These are, of course, the prime order subgroups of $G$.
In this paper we shall give an explicit formula for $x_G$ in terms
of the elements $x_{P_i}$, where $T=\{P_1,\cdots,P_m\}$ contains
exactly one subgroup from each conjugacy class of subgroups of $G$
of prime order. To state the formula we need the following
definitions: For $i=1,\cdots,m$, let $L_{P_i}$ be the set
$\{\{g_1,\cdots,g_{l_i}\}|G= \bigcup_{j=1}^{l_i}{g_jP_i}\}.$
The set $L_{P_i}$ has an obvious $G$-action, given by left multiplication.
The cartesian product $X=\prod_{i=1}^m{L_{P_i}}$ is a $G$-set via
the diagonal $G$-action. Let $\{w_j\}_{j=1,\cdots,s}$ be a set of
representatives of the different orbits of the action of $G$ on $X$.
Since $w_j\in X$ we can write
\begin{equation}\label{reps_formula}w_j=(\{g_{ij1},\ldots,g_{ijl_i}\})_{i=1}^m.\end{equation} The main result of this paper is the following:
\begin{proposition}\label{thm1} Suppose that,
for each $i=1,\cdots,m$, the element $x_{P_i}$ satisfies the equation
$tr_{P_i}(x_{P_i})=1$. Then the element
\begin{equation}\label{formula2}x_G =\sum_{j=1}^s{\Big[\prod_{i=1}^m{\Big(\prod_{t=1}^{l_i}g_{ijt}(x_{P_i})\Big)}\Big]}.\end{equation}
satisfies the equation $tr_G(x_G)=1$. \end{proposition} The $g_{ijt}$
in the proposition comes from formula \ref{reps_formula} above.
We will prove proposition \ref{thm1} in Section \ref{proof}.

As one might see Formula \ref{formula2} uses all the conjugacy
classes. One might ask if all the conjugacy classes are really
needed; that is, can we find a formula which uses only some of the
conjugacy classes? In a wider sense, one can ask if the existence or
non-existence of $x_P$, for $P$ in some of the conjugacy classes
affects the existence or non-existence of $x_P$ for $P$ in other
conjugacy classes. We will prove in this paper the following
proposition, which answers this question in a negative way:
\begin{proposition}\label{proposition2}
Let $T$ be a set of representatives of the
conjugacy classes of subgroups of prime order in $G$, and let $T = A
\amalg B$ be a partition of $T$ into two disjoint subsets, $A$ and $B$.
There exists a commutative ring $R = R_{A,B}$ upon which $G$ acts such that for every $P \in
A$ there is an element $x_P \in R$ with $tr_P(x_P) = 1$, and for
every $P \in B$ there exists no such $x_P$.\end{proposition}

The paper is arranged as follows. In the second section we shall
give some preliminaries which are needed for the rest of the paper.
In the third section we shall prove Proposition \ref{thm1}, and in the
fourth section we shall prove Proposition \ref{proposition2}.

\textbf{Acknowledgments.} The author would like to thank his
teachers Eli Aljadeff, Shlomo Gelaki, and Jack Sonn, for their
guidance, and for their great help in bringing this paper to its
present form. The author would also like to thank the referee
for his useful comments.
\end{section}

\begin{section}{PRELIMINARIES}\label{preliminaries}
Let $G$ be a finite group. A \textbf{$G$-ring} is a ring on which
$G$ acts by ring automorphisms. A \textbf{$G$-morphism} between two
$G$-rings $R$ and $S$ is a homomorphism of rings $\phi: R\rightarrow
S$ which is equivariant with respect to the $G$-action. A
\textbf{$G$-ideal} in a $G$-ring $R$ is an ideal $I\vartriangleleft
R$ such that for every $\sigma\in G$ and every $x\in I$,
$\sigma(x)\in I$. It is easy to see that if $\phi: R\rightarrow S$
is a $G$-morphism, then $ker\, \phi$ is a $G$-ideal, and if $I$ is a
$G$-ideal of the $G$-ring $R$, then the ring $R/I$ has a natural
$G$-ring structure such that the natural projection
$\pi:R\rightarrow R/I$ is a $G$-morphism.

If $A=\{a_1,\ldots,a_n\}$ is a subset of a $G$-ring $R$, then the
$G$-ideal generated by $A$ is the smallest $G$-ideal which contains
$A$. It is the same as the ideal generated by the elements $\{\sigma
a_i\}_{\sigma\in G, i=1,\ldots,n}$. We shall denote this ideal by
$(a_1,\ldots,a_n)_G$.

If $G$ acts on the set $X$, then the polynomial ring $\Z [X]$ is a
$G$-ring in a natural way- the action of $G$ on the indeterminates,
which are elements of $X$ is given. We extend this action uniquely
to an action of $G$ on the whole ring, using additivity and
multiplicativity.

We call $\Z[X]$ the \textbf{$G$-ring on the $G$-set $X$}. If $Y$ is any set,
then we can define $X=\{\sigma y\}_{\sigma\in G, y\in Y}$ with the
obvious $G$-action. The resulting $G$-ring $\Z [X]$ is called the
\textbf{free $G$-ring on $Y$}. Suppose that $X$ is a $G$-set,
$R$ is a $G$-ring, and that $\phi: X\rightarrow R$ is a map of
$G$-sets (i.e. equivariant with respect to the $G$-action). Then
there exists a unique $G$-morphism $\tilde{\phi}: \Z [X] \rightarrow
R$ which extends $\phi$. It is always true that there exists such a
unique ring homomorphism, and since $\phi$ is $G$-equivariant, it
follows that this ring homomorphism is a $G$-morphism.

We need a few basic facts about the trace map. Recall that if $H<G$
is a subgroup, then the trace map is defined as $$tr_H : k
\rightarrow k$$\begin{equation}tr_H(x) = \sum_{\sigma\in
H}{\sigma(x)}.\end{equation}
It is known that if $H<G$ and $tr_G$ is onto, then $tr_H$ is onto as well.
This can easily be seen by considering $x_H =
\sum_{i=1}^l{g_i(x_G)}$, where $g_1, \cdots, g_l$
are coset representatives of $H$ in $G$.
If $x_H$ is an element with trace 1 for the subgroup $H$, and
$g\in G$ is an arbitrary element, then $gx_H$ is an element with
trace 1 for the subgroup $gHg^{-1}$.
It follows that $H$ has an element with trace 1 if and only if
$gHg^{-1}$ has one. It is therefore suffices here to consider subgroups of $G$
only up to conjugacy.
If $|G|=p_1^{d_1}\ldots p_e^{d_e}$, then the existence of the
element $x_G$ is equivalent to the existence of the elements
$x_{P_i}$, where $P_i$ is a $p_i$-sylow subgroup of $G$ for
$i=1,\cdots,e$. This is because if $tr_G$ is onto, then $tr_{P_i}$
is onto, for every $i=1,\cdots,e$ as noted above, and if $x_{P_i}$
is an element with trace 1 for the group $P_i$, then
$tr_G(x_{P_i})=|G|/{p_i^{d_i}}$. Since these are coprime numbers, it
is easy to see that $tr_G$ is onto as well.
Another basic fact that will be needed in the sequel is this:
if $\phi:R \rightarrow S$ is a $G$-morphism between $G$-rings,
and $x_G$ is an element with trace 1 in $R$, then $\phi(x_G)$ is an
element with trace 1 in $S$. One can see this by considering the
trace of $\phi(x_G)$ in $S$ and using the fact that $\phi$ is
a $G$-morphism.
A question that arises naturally when one is searching a formula
for $x_G$ in terms of $x_H$, where $H$ varies over some set $T$ of
subgroups of $G$, is why would such a formula exists.
We need to know, of course, that the existence of $x_H$
for every $H$ in $T$ determines the existence of $x_G$.
In this case, the following proposition, which was proved in a
more general form by Shelah, is known. See \cite{Al2} and \cite{Al4} for a proof.
\begin{proposition} Let $G$ be a finite group and let $k$ be a
commutative ring upon which $G$ acts. Suppose that $T$ is a collection of subgroups of
$G$ such that the existence of $x_H$ for every $H\in T$ determines
the existence of $x_G$. Then there exist a polynomial formula for
$x_G$ in terms of the elements $\sigma x_H$, where $\sigma \in G$
and $H\in T$. \end{proposition}
\begin{remark} The formula is a universal
one, by which we mean that it does not depend on the particular ring $k$
and will work in any commutative ring. \end{remark}

We shall now explain why the surjectivity of the trace map is of
interest. Let $k$ be a $G$-ring. Denote the corresponding
homomorphism by $t:G\rightarrow\textrm{Aut}(k)$. Define the \emph{skew group
ring} $k_tG$ for this action to be the free $k$-module with
basis $\{u_\sigma\}_{\sigma\in G}$, whose multiplication is given by
the rule
$$a_{\sigma}u_{\sigma}a_{\tau}u_{\tau}=a_{\sigma}\sigma(a_{\tau})u_{\sigma\tau} \text{ for }\sigma,\tau\in G, \, a_{\sigma},a_{\tau}\in k.$$
Then $k$ has a natural $k_tG$-module structure given by the rule
$$\sum_{\sigma\in
G}{a_{\sigma}u_{\sigma}}\cdot{b}=\sum_{\sigma\in
G}{a_{\sigma}\sigma(b)}\text{ for }\sigma\in G, \, a_{\sigma},b\in k.$$ The following proposition explains the
connection between the structure of $k$ as a $k_tG$-module, and the
surjectivity of the map $tr_G$.
\begin{proposition}\label{proposition_projectivity} Let $k$, $G$, and $t$ be as above. Then $k$ is a projective
$k_tG$-module if and only if $tr_G$ is surjective. \end{proposition}
The proof of proposition \ref{proposition_projectivity} can be found in \cite{Al2}.
For a deeper treatment of skew group rings and modules over group rings,
see also \cite{Rim}, \cite{Cho}, and \cite{Sch}.

\end{section}

\begin{section}{PROOF OF PROPOSITION \ref{thm1}}\label{proof}
Recall the following notations from the introduction:
$T=\{P_1,\ldots,P_m\}$ is a set of representatives of conjugacy
classes of prime order
subgroups of $G$. $L_{P_i}=\{\{g_1,\cdots,g_l\}|G=
\bigcup_{i=1}^l{g_iH}\}$ is the set of sets of representatives of
the cosets of $P_i$ in $G$, upon which $G$ acts in the obvious way.
$X$ is the cartesian product $\prod_{i=1}^m{L_{P_i}}$ with the
diagonal $G$-action, and $w_j=(\{g_{ij1},\ldots,g_{ijl_i}\})_{i=1}^m
\textrm{ for } j=1,\cdots,s$ are representatives of the orbits of
the action of $G$ on $X$.

In the course of the proof of the formula, we will need the following
simple combinatorial principle:
\begin{equation}\label{Combinatorial}\prod_{i=1}^l{\Big(\sum_{j=1}^{t_i}{x_{ij}}\Big)} =
\sum_{j_1\leq t_1,\ldots,j_l\leq t_l}{\prod_{i=1}^l{x_{ij_i}}}.\end{equation}
This principle can easily be proved if one opens
the parenthesis in the left hand side of the equation.
We now state our main result:

\vspace{4pt}\textbf{Proposition \ref{thm1}} \textit{Suppose that,
for each $i=1,\cdots,m$, the element $x_{P_i}$ satisfies the equation
$tr_{P_i}(x_{P_i})=1$. Then the element
\begin{equation}\label{equ4}x_G =\sum_{j=1}^s{\Big[\prod_{i=1}^m{\Big(\prod_{t=1}^{l_i}g_{ijt}(x_{P_i})\Big)}\Big]}.\end{equation}
satisfies the equation $tr_G(x_G)=1$.}\vspace{4pt}

The rest of this section is devoted to the proof of Proposition
\ref{thm1}
\begin{proof}
For every $i$, we have the equation $tr_{P_i}(x_{P_i})=1$. This means
that \begin{equation}\label{eq1}\sum_{g \in
P_i}{g(x_{P_i})}=1.\end{equation} Acting with any $z \in G$ on the
last equation, we get \begin{equation}\label{eq2}\sum_{g \in
P_i}{zg(x_{P_i})}=1.\end{equation} Note that the last summation is
over all the elements in the coset $zP_i$. Since all these equations
equal $1$, their product also equals $1$. We first consider this product
for a fixed $i$. Using Equation \ref{Combinatorial}, we have
\begin{equation}\label{equ5}\prod_{z\in G/{P_i}}{\sum_{g\in P_i}{zg(x_{P_i})}}
= \sum_{\{g_1,\ldots,g_{l_i}\} \in
L_{P_i}}{\prod_{j=1}^{l_i}{g_j(x_{P_i})}}=1.\end{equation} The
meaning of $z\in G/{P_i}$ is that we took one element from each
coset of $P_i$ in $G$. It is easy to see that when we use Equation
\ref{Combinatorial} on the product in Equation \ref{equ5}, we get
summation over the different coset representatives of $P_i$ in $G$,
as indicated in the equation. Thus far, we have equation of the form
of Equation \ref{equ5} for each $i=1,\ldots,m$. Next, we shall
multiply all these equations together. Remember that
$X=\prod_{i=1}^m{L_{P_i}}$. Now by applying Equation
\ref{Combinatorial} again, we get:
\begin{equation}\label{equ6}\prod_{i=1}^m{\sum_{\{g_1,\ldots,g_{l_i}\} \in
L_{P_i}}{\prod_{j=1}^{l_i}{g_j(x_{P_i})}}}
=\sum_{(g_{i1},\ldots,g_{il_i})_{i=1}^m \in
X}{\prod_{i=1}^m{\prod_{j=1}^{l_i}{g_{ij}(x_{P_i})}}}=1.\end{equation}
Consider the set \begin{equation}A = \{
\prod_{i=1}^m{\prod_{j=1}^{l_i}{g_{ij}(x_{P_i})}}|
(g_{i1},\ldots,g_{il_i})_{i=1}^m \in X\}.\end{equation}
We consider
the elements of $A$ as formal products. The elements of $A$ can also
have an interpretation as elements of the commutative ring $k$. Since
$G$ acts on $X$, it also acts on $A$ by
\begin{equation}g\cdot\prod_{i=1}^m{\prod_{j=1}^{l_i}{g_{ij}(x_{P_i})}}=
\prod_{i=1}^m{\prod_{j=1}^{l_i}{gg_{ij}(x_{P_i})}}.\end{equation} We
claim the following:
\begin{lemma} For every $P_i \in T$, and every $a \in A$ we have
$stab(a) \cap P_i = 1$.\end{lemma}
\begin{proof}Let $1 \neq g \in P_i$, and let $a =
\prod_{i=1}^m{\prod_{j=1}^{l_i}{g_{ij}(x_{P_i})}}$. Consider the
element from the trivial coset of $P_i$ in $G$, $g_{i1}\in P_i$.
Call this the representative of the trivial coset of $P_i$ in $a$.
Since $g \in P_i$, the representative of the trivial coset of
$P_i$ in $g\cdot a$ is $gg_{i1}$. By assumption, $g\neq 1$, and
therefore $gg_{i1}\neq g_{i1}$. It follows that $g\cdot a \neq a$,
so $g \notin stab(a)$. Therefore $stab(a) \cap P_i = 1$ as
desired.
\end{proof}
Since stab($g\cdot a$) = $g\cdot stab(a)\cdot g^{-1}$, and the last lemma was
proved for an arbitrary element $a \in A$, we have the following
corollary:
\begin{corollary} For every subgroup $P<G$ of prime order, and
every $a \in A$, one has $stab(a) \cap P = 1$.\end{corollary}
\begin{proof}By the definition of $T$, there is an $i$ and a $g\in G$
such that $gPg^{-1}=P_i$. We have $g(stab(a)\cap P)g^{-1} =
gstab(a)g^{-1}\cap gPg^{-1} = stab(ga)\cap P_i=1$, and
therefore $stab(a)\cap P = 1$, as desired.\end{proof} Next, we
claim something stronger:
\begin{lemma} For every $a\in A, stab(a)=1$.\end{lemma}
\begin{proof} Suppose $stab(a)\neq 1$. Then $stab(a)$ contains a
subgroup of prime order $P$. The corollary above says that
$stab(a)\cap P=1$, and this is of course a contradiction.\end{proof}
We can now complete the proof of Proposition \ref{thm1}. We know that when we consider the
elements of $A$ as products in the ring $k$, we have
\begin{equation}\sum_{a \in A}a = 1\end{equation} (this is just Equation
\ref{equ6}). Let us choose a set of representatives of the different
orbits of the action of $G$ on $A$, which we shall denote by
$\{a_1,\ldots,a_q\}$. Since the stabilizer of each element in $A$ is
trivial, it follows that each element in $A$ can be written uniquely
in the form $\sigma a_j$, where $\sigma \in G, j=1,\ldots,q$. Thus
if we define $x_G = \sum_{j=1}^q{a_j}$, we have
\begin{equation}tr_G(x_G) = \sum_{\sigma \in G}{\sum_{j=1}^q{\sigma
a_j}} = \sum_{a \in A}{a} = 1.\end{equation} Thus $x_G$ is an
element with trace $1$ for the group $G$, and it is given by Formula \ref{equ4}, as desired.
\end{proof}
\begin{remark}
In order that the above formula will work, the ring $k$ need not be
commutative. It is enough that for every $P<G$ and every $\sigma \in
G$ we have $x_P\sigma(x_P) = \sigma(x_P)x_P$. Non-commutative rings
which satisfy this condition can be constructed artificially.
However, we do not know any natural examples of such rings. For the
general noncommutative case there is a formula in case the group $G$
is abelian, see \cite{Al3}.\end{remark}
\begin{remark} The only special property of the set $T$ we used in the course
of the proof above is the following: for every $1\neq g \in G$,
there exists a natural number $n$, an element $\sigma \in G$, and a
subgroup $P \in T$ such that $g^n \neq 1$ and $g^n \in \sigma P
\sigma^{-1}$. If we replace the set $T$ by any other set of subgroup
of $G$ which satisfies the above condition, then by the same proof we
will have a formula for $x_G$ in terms of the $x_N$'s where $N$ varies
over $T$.
\end{remark}
\end{section}

\begin{section}{INDEPENDENCE OF THE FACTORS}
In Section \ref{introduction} we gave a formula for $x_G$ in terms
of the elements $x_P$, where $P$ varies over the set
$T=\{P_1,\ldots,P_m\}$. The formula we gave uses \emph{all} the
$x_P$'s. One might ask if there exists a formula for $x_G$ which
does not use all the $x_P$'s. In a wider sense we can ask if it is
possible that the existence (or non-existence) of $x_{P_i}$ for some
of the subgroups $P_i\in T$ determines the existence (or
non-existence) of  some of the others. As we shall see here, this is
not the case. We will show here that the existence (or
non-existence) of some of the factors, does not say anything about
the existence (or non-existence) of the others. More precisely, we
shall prove the following result.

\vspace{4pt}\textbf{Proposition \ref{proposition2}} \textit{Let $T$
be as above. Suppose that $T$ is a disjoint union $T = A\amalg B$.
Then there exists a ring $R= R_{A,B}$ such that for every subgroup
$P\in A$ we have an element $x_P \in R$ with $tr_P(x_P)=1$, and for
every $P\in B$ there is not such an element.}\vspace{4pt}

Thus, we can view this proposition as saying that the existence of
the $x_P$'s for different $P's$ is independent. Moreover, this
proves that a formula with less factors cannot exist. Indeed, if we
would had a formula which uses only some of the groups in $T$, say,
only the subgroups in the subset $A \subset T$, then the existence
of $x_P$ for $P \in A$ would have implied the existence of $x_G$,
and this in turn would have implied the existence of $x_Q$ for every
$Q \in T\backslash A$, contradicting the above proposition.
\begin{proof}
We shall construct the ring $R_{A,B}$ explicitly. For every subgroup
$P \in B$ we have the $G$-set of the left cosets of $P$ in $G$. Let
us denote it by $Y_P = \{gP\}_{g \in G}$. Consider now the disjoint
union of these sets $Y = \amalg_{P \in B}{Y_P}$. Since $G$ acts on
each of the $Y_P$'s, $G$ acts on $Y$. Recall that the stabilizer for
this action is given by $stab(gP) = gPg^{-1}$. Now $Y$ is a $G$-set,
so we can build the $G$-ring on $Y$, $\mathbb{Z}[Y]$. Denote this
ring by $k$. Define $I = (\sum_{y \in Y}{y} - 1)$. Since the element
$\sum_{y \in Y}{y} - 1$ is $G$-invariant, $I$ is easily seen to be a
$G$-ideal, and we can consider the $G$-ring \begin{equation}R =
R_{A,B} = k/I.\end{equation} Let $P \in A$, and consider $R$ as a
$P$-ring. Since $G$ acts on the set $Y$, $P$ also acts on the set
$Y$. The stabilizer in $P$ of an arbitrary element of $Y$, $gQ$, is
$gQg^{-1} \cap P$. Since $T$ is a set of representatives of
conjugacy classes of subgroups of prime order, and $Q \in B$ , the
two subgroups $gQg^{-1}$ and $P$ are different, and therefore their
intersection is trivial. It follows that all the stabilizers for the
action of $P$ on $Y$ are trivial. Now choose a set of
representatives for the different orbits of the action of $P$ on
$Y$. Denote this set by $\{y_1,\ldots,y_e\}$. Define $y =
\sum_{i=1}^e{y_i}$. We thus have \begin{equation}tr_P(y) =
\sum_{\sigma \in P}{\sum_{i=1}^e{\sigma y_i}} = \sum_{y\in Y}y
\equiv 1 \ (mod \ I).\end{equation} The second equality follows from
the fact that $\{y_1,\ldots,y_e\}$ are representatives for the
action of $P$ on $Y$ and the stabilizers of this action are trivial.
The third equality follows from the definition of $I$. We see
therefore that $y+I \in R$ is an element of $R$ with $tr_P(y)=1$ (in
$R$). Now let $P \in B$. For every $\sigma \in P$ we have $\sigma P
= P$. Consider the map \begin{eqnarray}\phi : Y \rightarrow
\mathbb{Z}\\\nonumber a \mapsto \left\{ \begin{array}{ll}
            1 & \textrm{if $a = P$} \\
            0 & \textrm{else}
            \end{array} \right.\end{eqnarray}
Since $P$ (considered as an element of $Y$) is a fixed point for the
action of $P$ on $Y$, it can easily be seen that $\phi$ is a $P$-map
($\mathbb{Z}$ is considered to be a $P$-ring with the trivial action).
Therefore, it gives rise to a unique $P$-morphism $\tilde{\phi} :
\mathbb{Z}[Y] \rightarrow \mathbb{Z}$. We have
\begin{equation}\tilde{\phi}(\sum_{y\in Y}y - 1) = \sum_{y \in Y}\phi(y) - 1
= 1 + 0 + \ldots + 0 - 1 = 0.\end{equation} It follows that
$\tilde{\phi}(I)=0$ and therefore $\tilde{\phi}$ factors
through the natural projection $\pi : \mathbb{Z}[Y] \rightarrow R$.
This means that there exists a $P$-map $\psi : R \rightarrow
\mathbb{Z}$. Now suppose that there is an element $x_P \in R$ such
that $tr_P(x_P)=1$. Then $y_P = \psi(x_P) \in \mathbb{Z}$ satisfies
$tr_P(y_P)=1$. Since $\mathbb{Z}$ has the trivial $P$-action, and
$|P|$ is not invertible in $\mathbb{Z}$, such an element cannot
exist. In conclusion, for every $P\in A$ we have an element $x_P\in
R$ such that $tr_P(x_P)=1$, and for every element $P\in B$ we do not
have such an element. The ring $R = R_{A,B}$ therefore satisfies the
desired properties, and we are done. \end{proof}
\end{section}

\end{document}